\input amstex
\documentstyle{conm-p}
\NoBlackBoxes

\issueinfo{00}
  {}
  {}
  {XXXX}

\def\P{{\Bbb P}}

\def\E{{\Bbb E}}
\def\a{\alpha}
\def\d{\delta}

\def\modo#1{{\left|#1\right|}}
\def\normo#1{{\left\|#1\right\|}}
\def\smodo#1{{\mathopen|#1\mathclose|}}
\def\snormo#1{{\mathopen\|#1\mathclose\|}}

\def\striplenormo#1{\smodo{\smodo{\smodo{#1}}}}

\def\Proof:{\medskip\noindent{\bf Proof:}\ \ } \def\Proofof
#1:{\medskip\noindent{\bf Proof of #1:}\ \ }

\topmatter
\title A Note on Sums of Independent Random Variables\endtitle
\author Pawe{\l} Hitczenko and Stephen Montgomery-Smith\endauthor
\leftheadtext{PAWE{\L} HITCZENKO AND STEPHEN MONTGOMERY-SMITH}

\address Department of Mathematics,
North Carolina State University,
Raleigh, NC 27695--8205 \endaddress

\email pawel\@math.ncsu.edu\endemail

\urladdr http://math.ncsu.edu/\~{}pawel/\endurladdr

\thanks The first author was partially supported
by NSF grant DMS 9401345\endthanks

\address 
Department of Mathematics,
University of Missouri--Columbia, Columbia, MO 65211\endaddress

\email stephen\@math.missouri.edu\endemail

\urladdr http://math.missouri.edu/\~{}stephen\endurladdr

\thanks The second author was 
partially supported
by NSF grant DMS 9424396, and by the
University of Missouri Research Board
\endthanks



\abstract
In this note a two sided bound on the
tail probability of sums of independent, and either symmetric or
nonnegative, random variables is obtained. We utilize  a recent result by
Lata{\l}a
on bounds on moments of such sums. We also give a new proof of Lata{\l}a's
result for nonnegative random variables, and improve one of the constants
in his inequality.
\endabstract

\endtopmatter

\head 1. Introduction\endhead

Recently Lata\l a (1997) obtained the
following remarkable result: for a sequence of random variables $(X_n)$ and
$1\le p<\infty$ define the following
Orlicz norm
$$\striplenormo {(X_k)}_p=\inf\{\lambda>0: \prod_n\E|1+X_n/\lambda|^p\le
e^p\}.\leqno(1.1)$$ Lata\l a proved that $${e-1\over
2e^2}\striplenormo{(X_k)}_p\le \big(\E|\sum X_k|^p\big)^{1/p}\le
e\striplenormo{(X_k)}_p,\leqno(1.2)$$ provided $(X_n)$ are either symmetric
or positive, and in the first case $p\ge 2$, and in the second case $p\ge
1$. The main novelty here is the fact that, contrary to the classical
inequalities, the constants here are
independent of $p$. Certain particular cases of Lata\l a's result had been
known earlier (see e.g. Hitczenko
(1993), Gluskin and Kwapie\'n (1995) or Hitczenko, Montgomery-Smith and
Oleszkiewicz (1997)), but they
can be easily deduced from Lata{\l}a's inequality.

Of course, the ultimate goal is to obtain bounds on the tail probabilities
for sums of random variables. Lata{\l}a's result prompted us to investigate
that
problem. This program has been completed; our methods, which are based on
estimates for the decreasing
rearrangement of a random variable, work in a rather general setting. As a
result we were able to obtain
extensions of Lata{\l}a's result in various directions. The details of that
approach will be presented
elsewhere. The goal of this note is quite different; we will present a very
simple argument that allows
one to deduce tail bounds from Lata{\l}a's result. As a matter of fact,
this approach formally does not
really depend on Lata\l a's result, but it requires a knowledge of his
bounds on moments in order to be
employed successfully. We will also present a short proof (based on
decoupling techniques) of
Lata{\l}a's result for non-negative random variables. Our proof gives a
slightly better constant on the
left-hand side of (1.2).

Our notation is standard; for a sequence $(z_k)$ we let $\displaystyle
z_n^*=\max_{1\le k\le n}|z_k|$. The
letters $c$ and $C$ denote absolute constants whose values may change from
one use to the next.  We will write $S = \sum_{k=1}^\infty X_k$, and
$S_n = \sum_{k=1}^n X_k$, and $\normo S_p = (\E\modo S^p)^{1/p}$.

\head 2. Tail estimates via moment estimates\endhead

In this section we will to obtain two-sided estimates for tails of sums of
independent random variables. For the sake of brevity
we will concentrate on symmetric random variables, although it will be
clear that our arguments work for nonnegative random variables
 as well. In certain special cases tail inequalities have been obtained
from moment inequalities (see Gluskin and Kwapie\'n (1995), Hitczenko and
Kwapie\'n (1994)  or Hitczenko, Montgomery-Smith and Oleszkiewicz (1997)).
Also, in the case of multiples of Rademacher random variables, two-sided
estimates have been obtained by Montgomery and
Odlyzko (1988), and Montgomery-Smith (1990).

\proclaim{Theorem 2.1}
There exist positive constants
$c$, $C$, $\alpha $ and $\delta$
such that for all sequences of independent symmetric random variables
$(X_n)$, and for all $t$ such that
$$ t \ge {1\over 2} \snormo{\sum X_iI(|X_i|\le t)}_2
= {1\over 2} \left(\sum_{i=1}^n
\normo{X_iI(|X_i|\le t)}_2^2\right)^{1/2} ,$$ the following holds:
Let $p_t$ be the least $p$ such that
$$ \snormo{\sum X_iI(|X_i|\le t)}_p \ge 2 t .$$ Then we have the inequalities
$$ \P(|S_n| > t) \ge c\big\{\P(X_n^* > t)+\exp(-\a p_t)\big\},
\leqno(2.1)$$
and
$$ \P(|S_n| > 4 t)\le C\big\{\P(X_n^* > t)+\exp(-\delta
p_t)\big\}.\leqno(2.2) $$
If $t \le {1\over 2} \snormo{\sum X_iI(|X_i|\le t)}_2$, then $$ \P(|S_n| >
t) \ge c .$$
\endproclaim

\demo{Proof}
For a given $t$, let
$Y_i=X_iI(|X_i|\le t)$, and let
$s_n=\sum_{j=1}^n Y_j$.
Notice that $\normo{s_n}_p$ is a continuous, increasing function of $p$,
and that
$\normo{s_n}_2 \le 2t$. Hence either $2 \le p_t < \infty$ and
$ \snormo{s_n}_{p_t} = 2 t $,
or $p_t = \infty$ and $\normo{s_n}_\infty \le 2t$.

Let us start by proving $(2.1)$.
It follows from Levy's inequality and contraction principle
(see e.g. Kwapie\'n Woyczy\'nski (1992, Propositions 1.1.2 and 1.2.1)) that
$$ \P(X_n^*> t)\le 2\P(|S_n|> t),$$
and,
$$ \P(|s_n|> t)\le 2 \P(|S_n|> t).$$
Hence
$$\P(|S_n|>t)\ge {\textstyle{1\over4}}
\Big(\P(X_n^*>t)+\P(|s_n|>t)\Big).\leqno(2.3)$$ Now we can see that if $p_t
= \infty$, then the inequality is
established. In the case that $p_t < \infty$, we need to obtain a lower
estimate for the tail probability of a
maximum of partial sums of uniformly bounded symmetric random variables
$(Y_i)$. But for such random variables,
the following inequality is true (cf. Hitczenko (1994)): 
for all $q\ge p\ge1$, we have
$$\eqalign{\|s_n\|_q\le& C{q\over
p}\Big\{\|s_n\|_p+\|Y_n^*\|_q\Big\}\cr\le& C{q\over p}
\Big\{\|s_n\|_p+t\Big\}.\cr}\leqno(2.4)$$ We
also use the
Paley-Zygmund inequality that states that for any non-negative random
variable $Z$, and $0<\lambda<1$, $$ \P(Z > \lambda EZ) \ge (1-\lambda)^2
{(EZ)^2 \over EZ^2} .$$ Since $t = {1\over 2}\|s_n\|_{p_t}$ we have that $$
\eqalign{
\P(|s_n|> t)
&\ge \P(|s_n|^{p_t} > 2^{-p_t}\|s_n\|_{p_t}^{p_t})\cr &\ge
(1-2^{-p_t})^2{\|s_n\|_{p_t}^{2p_t}\over\|s_n\|_{2p_t}^{2p_t}}.\cr}$$ It
follows from (2.4) that the
denominator is no more than
$$C^{2p_t}\{\|s_n\|_{p_t}+t\}^{2p_t}\le
({\textstyle{3\over2}}C)^{2p_t}\|s_n\|_{p_t}^{2p_t}.$$ Therefore, we get the
estimate
$$\P(|s_n|\ge t)\ge
(1-2^{-p_t})^2({\textstyle{3\over2}}C)^ {-p_t}\ge \exp(-\a p_t),$$ which,
together with (2.3) gives (2.1).

Inequality (2.2) is an easy consequence of Chebyshev's inequality. If $p_t
< \infty$, then
$$\eqalign{\P(|S_n|>4t)\le& \P(X_n^*>t)+ \P(|s_n|>4t)\le
\P(X_n^*>t)+{E|s_n|^{p_t}\over (4t)^{p_t}}\cr\le&
\P(X_n^*>t)+2^{-p_t}=\P(X_n^*>t)+\exp(-\d p_t).\cr}$$ If $p_t = \infty$, we
use the same ideas, noticing that $\P(|s_n|>2t) = 0$.

Finally, if $t \le {1\over2} \normo{s_n}_2$, we apply the contraction
principle, the Paley-Zygmund
inequality, and $(2.4)$, to get
$$ \eqalignno{
2\P(|S_n| > t) &\ge \P(|s_n|^2 > {\textstyle{1\over4}} E|s_n|^2) \cr &\ge
{9\snormo{s_n}_2^4 \over 16\snormo{s_n}_4^4} \cr &\ge {9\snormo{s_n}_2^4
\over 16 C^4(\snormo{s_n}_2 + t)^4} \cr}$$ which is bounded below by a
universal constant.
\quad\qed
\enddemo

\remark{Remark} The above theorem allows us to approximate tails of the sums of
independent random
variables in terms of tails of the individual summands. This follows from
the fact that in view of
Lata\l a's result $p_t$ can be approximated using only information about
marginal distributions, and from the well known
inequality $${\sum \P(|X_i|> u)\over 1+\sum \P(|X_i|> u)}\le \P(X_n^*>
u)\le 2{\sum \P(|X_i|> u)\over 1+\sum \P(|X_i|> u)},\leqno(2.5)$$ which
gives tails of $X_n^*$ in terms of tails of individual summands.
\endremark

\head 3. Another proof of Lata{\l}a's result for nonnegative rv's\endhead

Here we intend to give another proof of Lata{\l}a's formula concerning
$\normo S_p$ for nonnegative random variables.

\proclaim{Theorem 3.1} Let $(X_n)$ be a sequence of positive independent
random variables. Then for all $p \ge 1$ we have that 
$$ \kappa \striplenormo{(X_n)}_p\le\snormo S_p
\le(e^p-1)^{1/p}
\striplenormo{(X_n)}_p ,\leqno(3.1)$$
where
$
\striplenormo{(X_n)}_p$ is given by (1.1), and $\kappa$ is the positive number
for which $f(\kappa) = e$, where
$$ f(x) = \sum_{k=0}^\infty {(2k+1)^k \over k!} x^k .$$
\endproclaim

\demo{Proof} First note that if $\striplenormo{(X_n)}_p \le 1$, then since $ 1 +
\sum_n X_n \le \prod_n (1+X_n)
$, we have that $\normo S_p^p \le e^p - 1$. This proves the second
inequality in (3.1) To prove the first,
we use certain results concerning decoupling. These ideas appear often in
the literature (usually in the context of mean-zero or symmetric random
variables, see e.g. Kwapie\'n and Woyczy\'nski (1992).) However, since we
will need
control of constants, we
cite the following, which is a special case of de la Pe\~na,
Montgomery-Smith and Szulga (1994, Theorem~2.1.)

\proclaim{Lemma~3.2} Let $(X_n)$ be a sequence of real valued independent
random variables. Let $(X_n^{(l)})$ be independent copies of $(X_n)$ for $1
\le l \le k$. Furthermore, let $f_{i_1,\dots,i_k}$ be elements of a Banach
space such that $f_{i_1,\dots,i_k} = 0$ unless the $i_1,\dots,i_k$ are
distinct. Then for any $1 \le p \le \infty$, we have that $$
\normo{\sum_{i_1,\dots,i_k} f_{i_1,\dots,i_k} X_{i_1} \cdots X_{i_k}}_p \le
(2k+1)^k \normo{\sum_{i_1,\dots,i_k} f_{i_1,\dots,i_k} X^{(1)}_{i_1} \cdots
X^{(k)}_{i_k}}_p .$$
\endproclaim

\noindent Now let us finish the proof of Theorem~3.1. Note that $$ \prod_n
\E\smodo{1+X_n}^p
= \normo{ \prod_n (1+X_n) }_p^p , $$
and so by Minkowski's inequality we have that $$ \normo{ \prod_n (1+X_n) }_p
\le
1+\sum_{k=1}^\infty
\normo{ \sum_{i_1 < \dots < i_k} X_{i_1} \cdots X_{i_k}}_p .$$ But if $k
\ge 1$ $$ \eqalignno{
\normo{ \sum_{i_1 < \dots < i_k} X_{i_1} \cdots X_{i_k}}_p &= {1\over k!}
\normo{\sum_{i_1 ,\dots, i_k\atop\hbox{\sevenrm distinct}} X_{i_1} \cdots
X_{i_k}}_p \cr
&\le
{(2k+1)^k\over k!} \normo{\sum_{i_1 ,\dots, i_k\atop\hbox{\sevenrm distinct}}
X_{i_1}^{(1)} \cdots X_{i_k}^{(k)}}_p \cr \noalign{\noindent (where
$(X_n^{(l)})$ are independent copies of $(X_n)$ for $1 \le l < \infty$)}
&\le
{(2k+1)^k \over k!} \normo{\sum_{i_1,\dots,i_k} X_{i_1}^{(1)} \cdots
X_{i_k}^{(k)}}_p \cr &= {(2k+1)^k \over k!} \normo S_p^k. \cr} $$ Hence
$$ \normo{ \prod_n (1+X_n) }_p
\le f(\normo S_p) ,$$
So, if
$ \normo S_p \le \kappa$,  then $$
\normo{ \prod_n (1+X_n) }_p \le e,$$ that is, $$ \striplenormo{(X_n)}_p \le
1.$$
\quad\qed
\enddemo

\remark{Remark} Our constant in the second inequality of (3.1) is
essentially the same as
Lata{\l}a's constant. But in the first inequality our constant,
which may numerically be shown to be about $0.1549$, is slightly
better than Lata{\l}a's constant, which is about
$0.1162$.
\endremark

\refstyle{1}


\widestnumber\key{11}
\Refs

\ref\key{1}
\by V.~de~la~Pe\~na, S.J.~Montgomery-Smith and J.~Szulga
\paper
Contraction
and decoupling inequalities for multilinear forms and U-statistics
\jour
Ann. of Probab.
\vol 22
\yr 1994
\pages 1745--1765
\endref

\ref\key{2}
\by E.D.~Gluskin and S.~Kwapie\'n
\paper
Tail and moment estimates for
sums of independent random variables
with logarithmically concave tails
\jour
Studia Math.
\vol 114
\yr 1995
\pages 303 - 309
\endref


\ref\key{3} 
\by P.~Hitczenko
\paper
Domination inequality for martingale transforms
of Rademacher sequence
\jour
Israel J. Math.
\vol 84
\yr 1993
\pages 161--178
\endref

\ref\key{4}
\by P.~Hitczenko
\paper
On a domination  of sums of random variables by
sums of conditionally independent ones
\jour Ann. Probab.
\vol 22
\yr 1994
\pages 453--468
\endref

\ref\key{5}
\by P.~Hitczenko and S.~Kwapie\'n
\paper On the Rademacher series
\inbook
Probability in Banach Spaces, Nine, Sandbjerg, Denmark
\ed J. Hoffmann-J\o rgensen, J. Kuelbs, M.B. Marcus
\publ Birkh\"auser
\publaddr Boston
\yr 1994
\pages 31--36
\endref

\ref\key{6} 
\by P.~Hitczenko, S.J.~Montgomery-Smith and K.~Oleszkiewicz
\paper Moment inequalities for sums of
certain independent symmetric random variables
\jour Studia Math
\vol 123
\yr 1997
\pages 15--42
\endref

\ref\key{7} 
\by S.~Kwapie\'n and W.A.~Woyczy\'nski
\book Random Series and Stochastic Integrals. Single and Multiple
\publ Birkh\"auser
\publaddr Boston
\yr 1992
\endref

\ref\key{8} 
\by R.~Lata\l a
\paper Estimation of moments of sums of independent
random variables
\jour Ann. Probab.
\vol 25
\yr 1997
\pages 1502--1513
\endref

\ref\key{9} 
\by H.L.~Montgomery and A.M.~Odlyzko
\paper Large deviations of sums of
independent random variables
\jour Acta Arithmetica
\vol 49
\yr 1988
\pages 427--434
\endref

\ref\key{10} 
\by S.J.~Montgomery-Smith
\paper The distribution of Rademacher sum
\jour Proc. Amer. Math. Soc.
\vol 109
\yr 1990
\pages 517--522
\endref

\endRefs
\enddocument